\documentclass[11pt]{article}
\input epsf
\usepackage{amsmath}
\usepackage{amssymb}
\usepackage{amscd}
\usepackage{epsfig}
\usepackage{amsthm}

\newtheorem{dfn}{Definition}[section]
\newtheorem{defn}[dfn]{Definition}

\newtheorem{rem}[dfn]{Remark}

\newtheorem{thm}[dfn]{Theorem}
\newtheorem{lem}[dfn]{Lemma}

\newtheorem{prop}[dfn]{Proposition}

\newtheorem{claim}[dfn]{Claim}
\newtheorem{cor}[dfn]{Corollary}

\newtheorem{cons}[dfn]{Construction}

\def\proof{\par\medskip\noindent{\it Proof: }}

\def\>{\rangle}
\def\<{\langle}

\def\3{\ss}
\def\8{\infty}

\begin{document}

% \size{12}{12pt}\selectfont

% \subjclass{57N10}

\title{Heegaard genus formula for Haken manifolds}
\author{Jennifer Schultens}

\maketitle
  
\begin{abstract}
Suppose $M$ is a compact orientable $3$-manifold and $Q \subset M$ a
properly embedded orientable boundary incompressible essential
surface.  Denote the completions of the components of $M - Q$ with
respect to the path metric by $M^1, \dots, M^k$.  Denote the smallest
possible genus of a Heegaard splitting of $M$, or $M^j$ respectively,
for which $\partial M$, or $\partial M^j$ respectively, is contained
in one compression body by $g(M, \partial M)$, or $g(M^j, \partial
M^j)$ respectively.  Denote the maximal number of non parallel
essential annuli that can be simultaneously embedded in $M^j$ by
$n_j$.  Then
\[g(M, \partial M) \geq \frac{1}{5}(\sum_j g(M^j, \partial M^j) - 
|M - Q| + 5 - 2\chi(\partial_-V) + 4\chi(Q) - 4\sum_j n_j)\] 
\end{abstract}

\maketitle
\vspace{2 mm}

Heegaard splittings have long been used in the study of $3$-manifolds.
One reason for their continued importance in this study is that the
Heegaard genus of a compact $3$-manifold has proven to capture the
topology of the $3$-manifold more accurately than many other
invariants.  In particular, it provides an upper bound for the rank of
the fundamental group of the $3$-manifold, and this upper bound need
not be sharp, as seen in the examples provided by M. Boileau and
H. Zieschang in \cite{BZ}.

We here prove the following: Let $M$ be a compact orientable
$3$-manifold and $Q \subset M$ an orientable boundary incompressible
essential surface.  Denote the completions of the components of $M -
Q$ with respect to the path metric by $M^1, \dots, M^k$.  Denote the
smallest possible genus of a Heegaard splitting of $M$, or $M^j$
respectively, for which $\partial M$, or $\partial M^j$ respectively,
is contained in one compression body by $g(M, \partial M)$ or $g(M^j,
\partial M^j)$ respectively.  Here $g(M, \partial M)$ is called the
\underline{relative genus} of $M$.  Denote the maximal number of non parallel
essential annuli that can be simultaneously embedded in $M^j$ by
$n_j$.  Then
\[g(M, \partial M) \geq \frac{1}{5}(\sum_j g(M^j, \partial M^j) - 
|M - Q| + 5 - 2\chi(\partial_-V) + 4\chi(Q) - 4\sum_j n_j)\] A
stronger inequality is obtained in the case in which $M$ and the
manifolds $M^j$ are acylindrical.

The formula derived in this paper provides a topological analogue to
the algebraic formula provided by R. Weidmann for the rank, i.e., the
minimal number of generators, of a group (see \cite{W}).  He proves
that if $G=A*\sb C B$ is a proper amalgamated product with malnormal
amalgam $C$, then \[{\rm rank}\,G\ge \frac{1}{3}({\rm rank}\,A+{\rm
rank}\,B-2\,{\rm rank}\,C+5)\]

The group $C < G$ is \underline{malnormal} if $gCg^{-1} \cap C =
\{1\}$ for all $g \; \epsilon \; G$.  Suppose that $M$ is a
$3$-manifold containing a separating incompressible surface $Q$ and
$Q$ cuts $M$ into two acylindrical $3$-manifolds $M^1 \sqcup M^2$.
Then Weidmann's formula tells us that the fundamental groups
$\pi_1(M), \pi_1(M^1), \pi_1(M^2), \pi_1(Q)$ satisfy the following
inequality:

\[r(\pi_1(M)) \geq \frac{1}{3}(r(\pi_1(M^1)) + r(\pi_1(M^2)) - 
2r(\pi_1(Q)) + 5)\]

The correspondence between the two results makes the formula derived
here particularly interesting, since it shows that the rank of the
fundamental group and the genus of a $3$-manifold satisfy a
similar linear inequality.  The construction and techniques used here
are a generalization of those in joint work with M. Scharlemann \cite{SS3}.
The complexity here is considerably more substantial.

In his book \cite{Jo}, K. Johannson derives a variant of the Heegaard
genus formula derived here for the special case in which $M$ is closed
and $M^j$ is acylindrical ({\it i.e.}, it contains no essential annuli or
tori) for $j = 1, \dots, k$:

\[\sum_j g(M^j) \leq 5g(M) + 2g(Q)\]

Which is equivalent to the following:

\[g(M) \geq \frac{1}{5}(\sum_j g(M^j) - 2 + \chi(Q))\]

See \cite[Proposition 23.40]{Jo}.  The inequality derived here applies
in greater generality.  In particular, Johannson's formula does not
apply to the interesting case of a surface bundle over the circle.

Section 2 of this paper shows how the generalized Heegaard splitting
of $M$ induces generalized Heegaard splittings of the submanifolds
$M^j$.  This construction provides the generalized Heegaard
splittings, but gives little control over their complexity.  Section 3
provides more specifics on the construction in Section 2 that provide
such control.  Section 4 proves the Main Theorem.  This proof
consists entirely of adding up and subtracting the appropriate numbers
from Sections 2 and 3.  

I wish to thank the MPIM-Bonn where this work was begun, RIMS, where
the final stage of this work was carried out and Professor Tsuyoshi
Kobayashi for inviting me to RIMS.  I would also like to thank
Professor Marty Scharlemann for helpful conversations.  This research
is supported in part by NSF grant DMS-9803826

\section{Preliminaries}

\vspace{2 mm}

For standard definitions concerning $3$-manifolds, see \cite{H} or
\cite{J}.

\begin{defn} 
For $L$ a properly embedded submanifold of $M$, we denote an open
regular neighborhood of $L$ in $M$ by $\eta(L)$ and a closed regular
neighborhood of $L$ by $N(L)$.  If $L$ is an orientable surface
properly embedded in $M$, then $N(L)$ is homeomorphic to $L \times
[-1, 1]$.  In this case we denote the subsets of $N(L)$ corresponding
to $L \times [-1, 0]$ and $L \times [0,-1]$ by $N_l(L)$ and $N_r(L)$
respectively.  We think of $L$ itself as corresponding to $L \times
\{0\}$.
\end{defn}

\begin{defn} \label{defn:cb} 
A \underline{compression body} is a $3$-manifold $W$ obtained from a
closed orientable surface $S$ by attaching $2$-handles to $S \times
\{-1\} \subset S \times I$ and capping off any resulting $2$-sphere
boundary components with $3$-handles.  We denote $S \times \{1\}$ by
$\partial_+W$ and $\partial W - \partial_+W$ by $\partial_-W$.
Dually, a compression body is an orientable $3$-manifold obtained from
a closed orientable surface $\partial_-W \times I$ or a $3-ball$ or a
union of the two by attaching $1$-handles.

In the case where $\partial_-W = \emptyset$ (i.e., in the case where a
$3-ball$ was used in the dual construction of $W$), we also call $W$ a
\underline{handlebody}.  In the case where $W = \partial_+W \times I$,
we call $W$ a \underline{trivial compression body}.

In the case of a compression body that is not connected, we further
require that all but one component of the compression body be a
trivial compression body.  The component of the compression body that
is non trivial is called the \underline{active component}.

In the following, we use the convention that $\chi(\emptyset) = 0$.
Define the \underline{index} of $W$ by $J(W) = \chi(\partial_-W) -
\chi(\partial_+W)$.  
\end{defn}

The index will usually be a positive integer, but the index of a
3-ball is -2 and the index of a solid torus or of a trivial
compression body is $0$.

\begin{defn}
A disk $D$ that is properly embedded in a compression body $W$ is
\underline{essential} if $\partial D$ is an essential curve in
$\partial_+W$.
\end{defn}

\begin{defn}
An annulus $A$ in a compression body $W$ is a \underline{spanning
annulus} if $A$ is isotopic to an annulus of the form
$(simple\;closed\;curve) \times I$ in the subset of $W$ homeomorphic
to $\partial_-W \times I$.
\end{defn}

\begin{defn} A \underline{set of defining disks} for a compression body
$W$ is a set of disks $\{D_1, \dots, D_n\}$ properly embedded in $W$
with $\partial D_i \subset \partial_+W$ for $i = 1$, $\dots, n$ such
that the result of cutting $W$ along $D_1 \cup \dots \cup D_n$ is
homeomorphic to $\partial_-W \times I$ along with a collection of
$3$-balls.
\end{defn}

\begin{defn} \label{defn:Heegaard splitting} A \underline{Heegaard
splitting} of a $3$-manifold $M$ is a decomposition $M = V \cup_S W$
in which $V$, $W$ are compression bodies such that $V \cap W =
\partial_+V = \partial_+W = S$.  We call $S$ the \underline{splitting
surface} or \underline{Heegaard surface}.

If $M$ is closed, the genus of $M$, denoted by $g(M)$, is the smallest
possible genus of the splitting surface of a Heegaard splitting for
$M$.  If $\partial M \neq \emptyset$, then the relative genus of $M$,
denoted by $g(M, \partial M)$, is the smallest possible genus of the
splitting surface of a Heegaard splitting for which $\partial M$ is
entirely contained in one of the compression bodies.
\end{defn}

\begin{defn} \label{defn:wr} 
A Heegaard splitting $M = V \cup_S W$ is \underline{reducible} if
there are essential disks $D_1 \subset V$ and $D_2 \subset W$, such
that $\partial D_1 = \partial D_2$.  A Heegaard splitting which is not
reducible is \underline{irreducible}.

A Heegaard splitting $M = V \cup_S W$ is \underline{weakly
reducible} if there are essential disks $D_1 \subset V$ and $D_2
\subset W$, such that $\partial
D_1 \cap \partial D_2 = \emptyset$.  A Heegaard splitting which is
not weakly reducible is \underline{strongly irreducible}.
\end{defn}

\begin{defn}
Let $M = V \cup_S W$ be a Heegaard splitting.  A Heegaard splitting is
\underline{stabilized} if there are disks $D \subset V, E \subset W$
such that $\vline \partial D \cap \partial E \vline = 1$.  A
\underline{destabilization} of $M = V \cup_S W$ is a Heegaard
splitting obtained from $M = V \cup_S W$ by cutting along the cocore
of a $1$-handle.  (E.g., if $D \subset V, E \subset W$ is a
stabilizing pair of disks, then $D$ is the cocore of a $1$-handle of
$V$ and the existence of $E$ guarantees that the result of cutting
along $D$ results in a Heegaard splitting.)  We say that a Heegaard
splitting $M = X \cup_T Y$ is a \underline{stabilization} of $M = V
\cup_S W$, if there is a sequence of Heegaard splittings $M = X^1
\cup_{T^1} Y^1, \dots, M = X^l \cup_{T^l} Y^l$ with $X^1 = X, Y^1 = Y,
X^l = V, Y^l = W$ and $X^r \cup_{T^r} Y^r$ is obtained from $X^{r-1}
\cup_{T^{r-1}} Y^{r-1}$ by a destabilization.
\end{defn}

The notion of strong irreducibility, due to Casson and Gordon in
\cite{CG}, prompted the following definition due to Scharlemann and
Thompson.

\begin{defn}
\label{general} A \underline{generalized Heegaard splitting} of a compact
orientable $3$-manifold $M$ is a decomposition $M = (V_1
\cup_{S_1} W_1) \cup_{F_1} (V_2 \cup_{S_2} W_2) \cup_{F_2} \dots
\cup_{F_{m-1}} (V_m \cup_{S_m} W_m)$.  Each of the $V_i$ and $W_i$ is
a compression body, $\partial_+V_i = S_i = \partial_+W_i$,
(i.e., $V_i \cup_{S_i} W_i$ is a Heegaard splitting of
a submanifold of $M$) and $\partial_-W_i = F_i = \partial_-V_{i+1}$.  We
say that a generalized Heegaard splitting is \underline{strongly
irreducible} if each Heegaard splitting $V_i \cup_{S_i} W_i$ is
strongly irreducible and each $F_i$ is incompressible in $M$.  We will
denote $\cup_i F_i$ by $\cal F$ and $\cup_i S_i$ by $\cal S$.  The
surfaces in $\cal F$ are called the \underline{thin levels} and the
surfaces in $\cal S$ the \underline{thick levels}.

Let $M = V \cup_S W$ be an irreducible Heegaard splitting.  We may
think of $M$ as being obtained from $\partial_-V \times I$ by
attaching all $1$-handles in $V$ followed by all $2$-handles in $W$,
followed, perhaps, by $3$-handles.  An \underline{untelescoping} of $M
= V \cup_S W$ is a rearrangement of the order in which the $1$-handles
of $V$ and the $2$-handles of $W$ are attached.  This rearrangement
yields a generalized Heegaard splitting.  If the untelescoping is
strongly irreducible, then it is called a \underline{weak reduction}
of $M = V \cup_S W$.  Here $\partial_-V_1 = \partial_-V$.  For
convenience, we will occasionally denote $\partial_-V_1$ by $F_0$ and
$\partial_-W_n$ by $F_n$.
\end{defn}

The Main Theorem in \cite{ST1} together with the calculation
\cite[Lemma 2]{SS2} implies the following:

\begin{thm} \label{thm:ind} Suppose $M$ is an irreducible compact
$3$-manifold.  Then $M$ possesses an unstabilized genus $g$ Heegaard
splitting \[M = V \cup_S W\] if and only if $M$ has a strongly
irreducible generalized Heegaard splitting
\[M = (V_1 \cup_{S_1} W_1) \cup_{F_1} (V_2 \cup_{S_2} W_2) \cup_{F_2}
\dots \cup_{F_{m-1}} (V_m \cup_{S_m} W_m)\] with $\partial_-V_1
= \partial_-V$ such that
\[\sum_{i=1}^m J(V_i) = 2g - 2 + 
\chi(\partial_-V).\]
\end{thm}

The details can be found in \cite{S1}.  Roughly speaking, one
implication comes from taking a weak reduction of a given Heegaard
splitting of genus $g$, the other from thinking of a given generalized
Heegaard splitting as a weak reduction of some Heegaard splitting.
The latter process is called the {\em amalgamation} (for details see
\cite{Sc1}) of the generalized Heegaard splitting.

A strongly irreducible Heegaard splitting can be isotoped so that its
splitting surface, $S$, intersects an incompressible surface, $P$, only in
curves essential in both $S$ and $P$.  This is a deep fact and is proven,
for instance, in  \cite[Lemma 6]{Sc2}.  This fact, together with the fact
that incompressible surfaces can be isotoped to meet only in essential
curves, establishes the following:

\begin{lem}  \label{lem:essential} Let $P$ be a properly embedded
incompressible surface in an irreducible $3$-manifold $M$ and let $M =
(V_1 \cup_{S_1} W_1) \cup_{F_1} \dots \cup_{F_{m-1}} (V_m \cup_{S_m}
W_m)$ be a strongly irreducible generalized Heegaard splitting of $M$.
Then ${\cal F} \cup {\cal S}$ can be isotoped to intersect $P$ only in
curves that are essential in both $P$ and ${\cal F} \cup {\cal S}$.
\end{lem}

\section{The Construction}

\vspace{2 mm} 

In this section, we suppose that $M$ is a compact orientable
$3$-manifold with generalized Heegaard splitting $M = (V_1 \cup_{S_1}
W_1) \cup_{F_1} \dots \cup_{F_{n-1}} (V_n \cup_{S_n} W_n)$ and $Q$ a
compact orientable surface in $M$ with $\partial Q \subset \partial M =
\partial_- V = \partial_-V_1$.  Assuming that $Q$ can be isotoped so
that all components of $Q \cap \cal (S \cup F)$ are essential in both
$Q$ and $\cal S \cup F$, we describe a construction for generalized
Heegaard splittings of the completion of the components of $M - Q$
with respect to the path metric.

In simpler contexts, a more concrete approach to this sort of
construction has been used.  A rough sketch is as follows: Consider a
Heegaard splitting $M = V \cup_S W$.  Make $V$ very thin, so it
intersects $Q$ in a collar of $\partial Q$ along with disks.  Then cut
along $Q$ and consider the completion of a component $C$ of $M - Q$.
Add a collar of the copies of components of $Q$ in the boundary of $C$
to $V \cap C$.  It is a nontrivial fact that this would indeed yield a
Heegaard splitting of $C$ and is the basis for the construction in
\cite{Jo}.  However, the approach there, though simpler at the outset,
requires far more work to gain control over the number
of components in $V \cap Q$.  Thus, although the construction in this
section appears to be more complicated than required, it will become
evident in the following sections that it allows for more satisfactory
control over the sum of indices of the resulting generalized Heegaard
splittings.

In the following, we will abuse notation slightly and consider $Q
\times [-1,1]$ to be lying in $M$ with $Q = Q \times \{0\}$ via the
homeomorphisn with $N(Q)$.  We will further assume that $F_i \cap (Q
\times [-1,1]) = (F_i \cap Q) \times [-1,1]$ and similarly for $S_i$,
for all $i$.  Denote the completions of the components of $M - Q$ with
respect to the path metric, {\it i.e.}, the $3$-manifolds into which
$Q$ cuts $M$, by $M^1, \dots, M^k$.  Note that $M^1 \sqcup \dots
\sqcup M^k$ is homeomorphic to $M - (Q \times (-1, 1))$.

\begin{defn}
A properly embedded surface $Q$ in a $3$-manifold $M$ is
\underline{essential} if it is incompressible
and not boundary parallel.
\end{defn}

\begin{rem}
An essential surface can be boundary compressible.  Recall that if a
surface $Q$ in a $3$-manifold $M$ is boundary compressible, then there
is a disk $D$ in $M$ such that $interior(D) \cap Q = \emptyset$ and
such that boundary $\partial D = a \cup b$ with $a, b$ connected arcs
and $a \subset Q$ and $b \subset \partial M$.  Supposing that $Q$ is
boundary compressible in $M$, then $D$ provides instructions for
modifying $Q$.  Specifically, replace a small collar of $a$ in $Q$ by
two parallel copies of $D$.  This modification is called a
\underline{boundary compression} of $Q$ along $D$.  Here $D$ is called
a \underline{boundary compressing disk} for $Q$.
\end{rem}

\begin{defn}
The two copies of $Q$ in $\partial (M^1 \sqcup \dots \sqcup M^k)$ are
called the \underline{remnants of Q}.
\end{defn}

\begin{defn}
Let $F$ be a closed orientable surface.  A generalized compression
body is an orientable 3-manifold $W$ obtained from $F \times I$ or a
3-ball or a union of the two by attaching 1-handles.  If attached to
$F \times I$, the 1-handles must be attached to $F \times \{1\}$.

We denote $F \times \{-1\}$ by $\partial_-W$.  We denote $\partial F
\times I$ by $\partial_v W$ and $\partial W - (\partial_-W \cup
\partial_v W)$ by $\partial_+W$.

A set of defining disks for $W$ is a set of disks ${\cal D}$ with
boundary in $\partial_+W$ that cut $W$ into $\partial_-W \times I$
together with a collection of $3$-balls.
\end{defn}

\begin{lem} \label{lem:auxcut}
Suppose that $W$ is a generalized compression body and $Q \subset W$
is a properly embedded connected incompressible surface disjoint from
$\partial_v W$.  Suppose further that $Q$ meets $\partial_+W$ and that
$\chi(Q) \leq 0$. Then either $Q$ is a spanning annulus, or
there is a boundary compressing disk $D$ for $Q$ such that $\partial D
\cap \partial W \subset \partial_+W$.
\end{lem}

\proof Let ${\cal D}$ be a set of defining disks for $W$.  Since $Q$
is incompressible, an innermost disk argument shows that $Q$ can be
isotoped so that it intersects the components of ${\cal D}$ in arcs.
Furthermore, an outermost arc argument shows that, after isotopy, any
such arc of intersection is essential in $Q$.  Now if there are arcs
of intersection, then we choose one that is outermost in ${\cal D}$
and see that the outermost disk it cuts off is a boundary compressing
disk for $Q$.

If there are no such arcs of intersection, then we cut along ${\cal
D}$ to obtain a 3-manifold homeomorphic to $\partial_-W \times I$.  It
is well known that an incompressible and boundary incompressible
surface in a product is either horizontal or vertical.
Here the horizontal case is ruled out because $Q$ does not meet
$\partial_v W$.  Thus $Q$ is either a spanning annulus or
boundary compressible.

Now suppose that $Q$ is boundary compressible via a disk $\tilde D$
such that $\partial \tilde D \cap \partial W \subset \partial_-W$.
Then since $\partial Q$ meets $\partial_+W$ there is a dual boundary
compressing disk $D$ such that $\partial D \cap \partial W \subset
\partial_+W$, as required.  \qed

\begin{rem} \label{rem:RS}
The construction here is relevant in the case in which for each component
$Q_c$ of $Q$, $Q_c \cap({\cal F \cup S}) \neq \emptyset$.  If, on the
other hand, there is a component $Q_c$ of $Q$ for which $Q_c
\cap({\cal F \cup S}) = \emptyset$, then we may treat this component
separately.  As it lies entirely in one of the compression bodies
$V_1, W_1, \dots, V_n, W_n$, it must in fact be parallel to a
component of ${\cal F}$.  If components of $Q$ are parallel into
components of $\cal F$, then a much simpler construction yields a
stronger result, see Proposition \ref{prop:par}.
\end{rem}

\begin{lem} \label{lem:cut}
Suppose $M = (V_1 \cup_{S_1} W_1 ) \cup_{F_1} \dots \cup_{F_{n-1}}
(V_n \cup_{S_n} W_n)$ is a generalized Heegaard splitting and suppose
$Q \subset M$ is an essential boundary incompressible surface.  Also
suppose that no component of $Q$ is parallel into ${\cal F}$.  Suppose
furthermore that $Q$ has been isotoped so that all components of $Q
\cap ({\cal S} \cup {\cal F})$ are essential in both $Q$ and ${\cal S}
\cup {\cal F}$ and so that the number of components in this
intersection is minimal subject to this condition.  Then for each $i$,
each component of the completion of $V_i - Q$ and $W_i - Q$ with
respect to the path metric is a generalized compression body.
\end{lem}

\proof
Note that under the above assumptions there will be no component of $Q
\cap V_i$ that does not meet $\partial_+V_i$, for such a component
would be parallel into $\partial_-V_i$.  Thus each component of $Q
\cap V_i$ satisfies the hypotheses of Lemma \ref{lem:auxcut}.
Similarly for $Q \cap W_i$.

Let $\tilde Q$ be a component of $Q \cap V_i$.  Since each component
of $\tilde Q \cap \partial_+V_i \subset Q \cap \partial_+V_i$ is
essential in $\partial_+V_i$ and in $Q$, $\tilde Q$ is not a disk.  If
$\tilde Q$ is a spanning annulus, then we may cut along this spanning
annulus and obtain a generalized compression body.  If $\tilde Q$ is
boundary compressible via a boundary compressing disk that meets
$\partial_+V_i$, then we may perform the boundary compression along
this disk to obtain $\tilde Q_b$.  The components of the completion of
$V_i - \tilde Q$ with respect to the path metric can be obtained from
the components of the completion of $V_i - \tilde Q_b$ with respect to
the path metric by attaching a 1-handle with cocore the boundary
compressing disk.

We prove the lemma by induction on $-\chi(Q \cap V_i)$.  This is
accomplished by repeated application of the argument above.  The same
holds for $Q \cap W_i$.  \qed

\begin{defn}
For a submanifold $N \subset M$, we will denote $N \cap M^j$ by $N^j$.
E.g., $S_3^1 = S_3 \cap M^1$, $W_5^2 = W_5 \cap M^2$, $Q^j = Q \cap
M^j$.  
\end{defn}

Note that $S_i^j$ and $F_i^j$ will typically not be closed surfaces.
Also, $V_i^j$ and $W_i^j$ will typically not be compression bodies,
only generalized compression bodies.  The following construction
appears to be a fairly natural way of ``capping off'' the components
of $\partial_v V_i^j$ and $\partial_v W_i^j$ with appropriate
((punctured surface) $\times I)$'s.  This is the first step in
constructing generalized Heegaard splittings on the submanifolds
$M^j$.  The difficulty lies in ``capping off'' $\partial_v V_i^j$ and
$\partial_v W_i^j$ in a way that is consistent.

\vspace{2 mm}

\begin{cons} \label{cons:key} 
(The Main Construction) Let $M$ be a compact possibly closed
orientable irreducible $3$-manifold.  Let \[M = (V_1 \cup_{S_1} W_1)
\cup_{F_1} \dots \cup_{F_{n-1}} (V_n \cup_{S_n} W_n)\] be a strongly
irreducible generalized Heegaard splitting of $M$.  Let $Q$ be a
compact possible closed orientable essential boundary incompressible
not necessarily connected surface properly embedded in $M$.  Denote
$\partial_-V_1$ by $\partial_-M$ and $\partial_-W_n$ by $\partial_+M$.
Then $\partial M = \partial_-M \cup \partial_+M$.  Suppose that
$\partial Q \subset \partial_-M$.  Suppose further that no component
of $Q$ is parallel to a component of $\cal F$ and that $\cal S \cup
\cal F$ has been isotoped so that all components of $Q \cap (\cal S
\cup \cal F)$ are essential in both $Q$ and $\cal S \cup \cal F$ and
so that the number of such components of intersection is minimal.

Denote the completions of the components of $M - Q$ with respect to
the path metric by $M^1, \dots, M^k$.  We construct generalized
Heegaard splittings for $M^1, \dots, M^k$, respectively, from \[M =
(V_1 \cup_{S_1} W_1) \cup_{F_1} \dots \cup_{F_{n-1}} (V_n \cup_{S_n}
W_n)\] We call these generalized Heegaard splittings the induced
Heegaard splittings of $M^1, \dots, M^k$, respectively.
\end{cons}

\begin{figure}
\centerline{\epsfxsize=2.0in \epsfbox{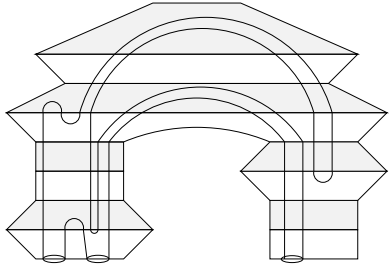}}
\caption{The surface $Q$ in the generalized Heegaard splitting}
\label{H0}
\end{figure}

Let $h:M \rightarrow [0,1]$ be a Morse function corresponding to $M =
(V_1 \cup_{S_1} W_1) \cup_{F_1} \dots \cup_{F_{n-1}} (V_n \cup_{S_n}
W_n)$ with $S_i = h^{-1}(s_i)$ and $F_i = h^{-1}(f_i)$ for appropriate
$s_1, \dots, s_n$ and $f_1, \dots, f_{n-1}$.  Note that $s_1, \dots,
s_n, f_1, \dots, f_{n-1}$ are regular values of $h$.  Our assumptions
on ${\cal S \cup F}$ guarantee that $s_1, \dots, s_n, f_1, \dots,
f_{n-1}$ are also regular values of $h|_Q$.

Set $q_t = Q \cap h^{-1}(t)$ and $Q_{[t_1,t_2]} = Q \cap
h^{-1}([t_1,t_2])$.  Here $q_t$ will consist of a collection of
circles (at least for regular values of $h|_Q$) and $Q_{[t_1, t_2]}$
will be a subsurface of $Q$ with $\partial Q_{[t_1, t_2]} = q_{t_1}
\cup q_{t_2}$.  Consider a bicollar $Q \times [-1, 1]$ of $Q$ in $M$.
If $q \in (Q \cap {\cal S \cup F})$, then we assume that $h(q,t) =
h(q)$ for all $t \in I$.

For all $i$, set \[Q_{s_i}^- = (q_{s_i} \times [-1, -1 +
\frac{1}{2i}]) \cup (Q_{[s_i, 1]} \times \{-1 +
\frac{1}{2i}\}) \] \[Q_{s_i}^+ = (q_{s_i} \times [1 - \frac{1}{2i}, 1]) 
\cup (Q_{[s_i, 1]} \times \{1 - \frac{1}{2i}\}) \]

For all $0 < i < n$, set \[Q_{f_i}^- = (q_{f_i} \times [-1, -1 +
\frac{1}{2i+1}]) \cup (Q_{[f_i, 1]} \times \{-1 + \frac{1}{2i+1}\})
\]  \[Q_{f_i}^+ = (q_{f_i} \times [1 - \frac{1}{2i+1}, 1]) \cup (Q_{[f_i, 1]} 
\times \{1 - \frac{1}{2i+1}\}) \] Recall that $F_0 = \partial_-V_1$
and $F_n = \partial_-W_n$.  Set \[Q_{f_0}^- = (q_{f_0} \times [-1,
-\frac{1}{100}]) \cup (Q \times \{-\frac{1}{100}\}) \] \[Q_{f_0}^+ =
(q_{f_0} \times [\frac{1}{100}, 1]) \cup (Q \times
\{\frac{1}{100}\}) \] \[Q_{f_n}^{\pm} = \emptyset\]

For all $i$, set \[\tilde F_i = (F_i - (F_i \cap (Q \times [-1, 1])))
\cup Q_{f_i}^+ \cup Q_{f_i}^-\] and \[\tilde S_i = (S_i - (S_i \cap (Q
\times [-1, 1]))) \cup Q_{s_i}^+ \cup Q_{s_i}^-\] Then for all $i$,
$\tilde F_i$ and $\tilde S_i$ are closed surfaces.  Note in particular
that since $F_n$ does not meet $Q$, we have $\tilde F_n = F_n$.  Let
$\tilde V_i$ be the cobordism between $\tilde F_{i-1}$ and $\tilde
S_i$ for $i = 1, \dots, n$ and let $\tilde W_i$ be the cobordism
between $\tilde S_i$ and $\tilde F_i$ for $i = 1, \dots, n$. Here
neither $\tilde V_i^j$ nor $\tilde W_i^j$ need be a compression body.

\begin{figure}
\centerline{\epsfxsize=2.0in \epsfbox{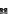}}
\caption{A collection of arcs $\alpha_i$ that cuts $Q_{[f_{i-1},
s_i]}$ into a spanning annulus}
\label{alphas}
\end{figure}

Let $\alpha_i$ be a union of properly embedded arcs in
$Q_{[f_{i-1},s_i]}$ disjoint from $q_{f_{i-1}}$ that cut
$Q_{[f_{i-1},s_i]}$ into disks and spanning annuli, see Figure
\ref{alphas}.  I.e., $Q_{[f_{i-1},s_i]} - \alpha_i$ is homeomorphic to
$(q_{f_{i-1}} \times I) \cup$ (disks).  Analogously, choose $\beta_i$
in $Q_{[s_i,f_i]}$.  Do this in such a way that $\partial \alpha_i
\cap \partial \beta_i = \emptyset$.  Then, set $\dot V_i^j = (\tilde
V_i^j - (\eta(\alpha_i \times \{\pm 1\})^j) \cup (\eta(\beta_i \times
\{\pm 1\})^j)$ and $\dot W_i^j = (\tilde W_i^j - (\eta(\beta_i \times
\{\pm 1\})^j) \cup (\eta(\alpha_i \times \{\pm 1\})^j)$.

\begin{claim} $\dot V_i^j$ and $\dot W_i^j$ are compression 
bodies for $i = 1, \dots, n$ and $j = 1, \dots, k$.
\end{claim}

\vspace{1 mm} Consider the construction of $\dot V_i^j$.  Since each
component of $Q \cap (\cal S \cup \cal F)$ is essential in both $Q$
and $\cal S \cup \cal F$ and since the number of such components is
minimal, $Q_{[f_{i-1}, s_i]} = Q \cap V_i$ is an essential surface.
Note that $Q_{[f_{i-1}, s_i]}$ might be boundary compressible, but
that the minimality assumption on the number of components of $Q \cap
(\cal S \cup \cal F)$ guarantees that the hypotheses of Lemma
\ref{lem:cut} are met.  Hence $V_i^j$ is a generalized compression
body.  Here $\tilde V_i^j$ is obtained from $V_i^j$ by adjoining a
3-manifold homeomorphic to $Q_{[f_{i-1}, 1]} \times I$ along a
subsurface homeomorphic to $Q_{[f_{i-1}, s_i]}$.  Furthermore $\dot
V_i^j$ is obtained from $V_i^j$ by adjoining the same 3-manifold but
along $\partial_v V_i^j = (q_{f_i} \times \{\pm 1\}) \times I$ and
disks and attaching 1-handles.  The result is thus a compression body.
Similarly for $\dot W_i^j$.

\vspace{1 mm} Set $\dot S_i^j = \partial_+ \dot V_i^j = \partial_+
\dot W_i^j$ and $\dot F_i^j = \partial_- \dot W_i^j = \partial_-\dot
V_{i+1}^j$ for all $i, j$.  Further set $\dot F_0^j = \partial_- \dot
V_1^j$ and $\dot F_n^j = \partial_- \dot W_n^j$ for all $j$.  The
generalized Heegaard splitting induced on $M^j$ is $M^j = (\dot V_1^j
\cup_{\dot S_1^j} \dot W_1^j) \cup_{\dot F_1^j} \dots \cup_{\dot
F_{n-1}^j} (\dot V_n^j \cup_{\dot S_n^j} \dot W_n^j)$.

Strictly speaking, the compression bodies may have to be relabelled.
For recall that if a compression body is not connected, then it has
exactly one active component.  Suppose, for instance, that $\dot
V_1^1$ is not connected and that more than one component is non
trivial.  (Note that it follows that $\dot W_1^1$ is also not connected.  For
$|\dot V_1^1| = |\partial_+ \dot V_1^1| = |\partial_+ \dot W_1^1| =
|\dot W_1^1|$.)  We then insert trivial compression bodies and relabel
as necessary.

\begin{rem}
Despite the possible insertion of trivial compression bodies and
relabelling, we maintain the notation $M^j = (\dot V_1^j \cup_{\dot
S_1^j} \dot W_1^j) \cup_{\dot F_1^j} \dots \cup_{\dot F_{n-1}^j} (\dot
V_n^j \cup_{\dot S_n^j} \dot W_n^j)$.  This helps us to keep track of
the relation to the original generalized Heegaard splitting on $M$.
This is especially helpful in the computations below.
\end{rem}

\begin{lem} \label{lem:construction} (Preliminary Count)
For the construction above, \[\sum_j\sum_i J(\dot V_i^j) =
\sum_i(J(V_i) + 2\chi(Q_{[f_{i-1}, s_i]}) + 4(|\alpha_i| +
|\beta_i|))\]
\end{lem}

\proof
Consider the Euler characteristics of the surfaces in the construction
above.  Since the Euler characteristic of a circle is $0$,
$\sum_j\chi(S_i^j) = \chi(S_i)$, for $i = 1, \dots, n$ and
$\sum_j\chi(F_i^j) = \chi(F_i)$, for $i = 0, \dots, n-1$.
Furthermore, $\sum_j \chi(\tilde S_i^j) = \sum_j \chi(S_i^j) +
2\chi(Q_{[s_i, 1]})$ and $\sum_j\chi(\tilde F_i^j) = \sum_j\chi(F_i^j)
+ 2\chi(Q_{[f_i, 1]})$.  Thus $\sum_j \chi(\dot S_i^j) = \sum_j
\chi(S_i^j) + 2\chi(Q_{[s_i, 1]}) - 4|\alpha_i| - 4|\beta_i|$, and
$\sum_j \chi(\dot F_i^j) = \sum_j \chi(F_i^j) + 2\chi(Q_{[f_i, 1]})$.
Hence
\[\sum_i \sum_j J(\dot V_i^j) = \]
\[\sum_i \sum_j (\chi(\dot F_{i- 1}^j) - \chi(\dot
S_i^j)) = \]
\[\sum_i \sum_j (\chi(F_{i-1}^j) - \chi(S_i^j)) + 2(
\chi(Q_{[f_{i-1},1]}) - \chi(Q_{[s_i, 1]}) + 4(|\alpha_i| + 
|\beta_i|)) = \]
\[\sum_i (\chi(F_{i-1}) - \chi(S_i) + 
2 \chi(Q_{[f_{i-1},s_i]}) + 4(|\alpha_i| + |\beta_i|)) = \]
\[\sum_i (J(V_i) + 2\chi(Q_{[f_{i-1},s_i]}) +
4(|\alpha_i| + |\beta_i|))\]
\qed

\section{How many arcs do we need?}

\vspace{2 mm}

In order to perform the required calculations we must count the number
of arcs required for $\alpha_i, \beta_i$ in the Main Construction.  To
do so, we define $\alpha_i, \beta_i$ more systematically.  

\begin{lem} \label{lem:arcs1}
In the Main Construction we may choose $\alpha_i$ so
that the number of components of $\alpha_i$ is equal to
$-\chi(Q_{[f_{i-1}, s_i]}) + d_i$, where $d_i$ is the number of
components of $Q_{[f_{i-1}, s_i]}$ that do not meet $F_{i-1}$.
\end{lem}

\proof
Recall that we are assuming that $\cal S \cup \cal F$ has been
isotoped so that all components of $Q \cap (\cal S \cup \cal F)$ are
essential in both $Q$ and $\cal S \cup \cal F$ and so that the number
of such components of intersection is minimal.

We proceed by induction on $-\chi(Q_{[f_{i-1}, s_i]})$.  If
$-\chi(Q_{[f_{i-1}, s_i]}) = 0$, then $Q_{[f_{i-1}, s_i]}$ consists of
annuli.  If a component $\tilde A$ of $Q_{[f_{i-1}, s_i]}$ is not a
spanning annulus, then the above minimality assumption implies that
its boundary lies in $\partial_+V_i$.  Thus Lemma \ref{lem:auxcut}
locates a boundary compressing disk $D$ for $\tilde A$ such that
$\partial D = a \cup b$ with $a \subset \tilde A$ and $b \subset
\partial_+V_i$.  In this case $\alpha_i$ consists of the arcs $a$ for
each such annulus.  Take $\alpha^0$ to be the set of such spanning
arcs.

Suppose $-\chi(Q_{[f_{i-1}, s_i]}) > 0$ and let $\tilde Q$ be a
component of $Q_{[f_{i-1}, s_i]}$ for which $-\chi(\tilde Q) > 0$.
Again, Lemma \ref{lem:auxcut} locates a boundary compressing disk $D$
for $\tilde Q$ such that $\partial D = a \cup b$ with $a \subset
\tilde Q$ and $b \subset \partial_+V_i$.  Set $\alpha^n = a$.

Consider the surface $Q^1$ obtained from $Q_{[f_{i-1}, s_i]}$ by a
boundary compression along $D$.  Then $-\chi(Q^1) = -\chi(Q_{[f_{i-1},
s_i]}) - 1$.  By inductive hypothesis there is a collection of arcs
$\alpha^0, \dots, \alpha^{n-1}$ such that the complement of $\alpha^0
\cup \dots \cup \alpha^{n-1}$ in $Q^1$ consists of spanning annuli and
disks.  Now take $\alpha_i$ to be the collection $\alpha^0, \dots,
\alpha^n$.

The fact that the number of components of $\alpha_i$ is equal to
$-\chi(Q_{[f_{i-1}, s_i]}) + d_i$ follows from our choice of
$\alpha_i$ as a subcollection of the arcs in $\alpha_1$ that cuts
$Q_{[f_{i-1}, s_i]}$ into spanning annuli and disks and that has the
minimal number of components among all such subcollections.
\qed

\begin{figure}
\centerline{\epsfxsize=3.0in \epsfbox{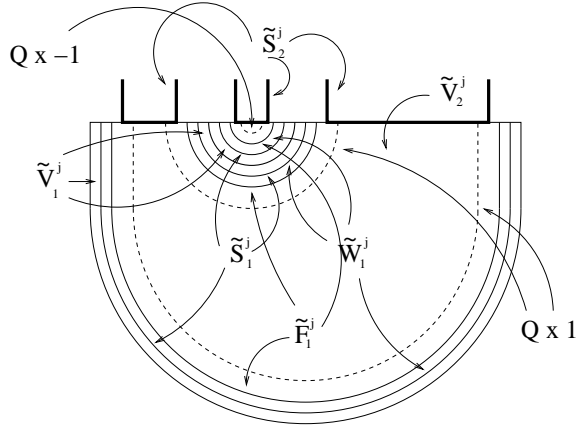}}
\caption{The construction of $\tilde V_2^j$}
\label{updown}
\end{figure}

\begin{figure}
\centerline{\epsfxsize=3.0in \epsfbox{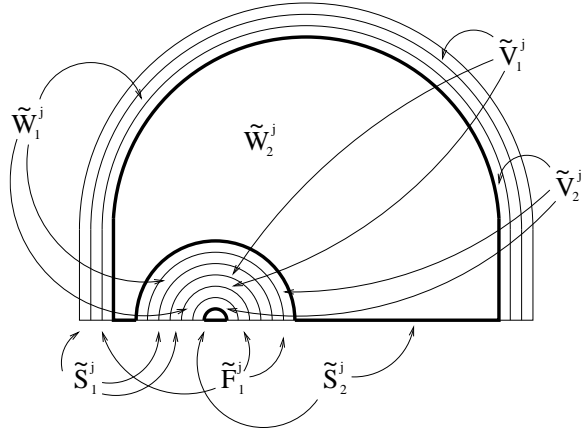}}
\caption{The construction of $\tilde W_2^j$}
\label{up}
\end{figure}

The same strategy could be used to locate a collection of arcs
$\beta_i$.  But there is a crucial difference between
$Q_{[f_{i-1}, s_i]}$ and $Q_{[s_i, f_i]}$.  This asymmetry in our
construction may be exploited to show that in fact all arcs in
$\beta_i$ are superfluous.  See Figures \ref{updown} and \ref{up}.

\begin{lem} \label{lem:arcs2}
In the Main Construction, we may choose, $\beta_i = \emptyset$.
\end{lem}

\proof Consider a collection of arcs $\beta_i$ in $Q_{[s_i, f_i]}$
constructed via the argument in Lemma \ref{lem:arcs1}.  Observe that
each arc found in Lemma \ref{lem:arcs1} was part of the boundary of a
boundary compressing disk $\tilde D$.  Here $\partial \tilde D = b
\cup c$ with $b \in Q_{[s_i, f_i]}$ and $c \in \partial_+W_i$. Denote
the collection of boundary compressing disks corresponding to the
components of $\beta_i$ by ${\cal D}_{\beta_i}$.  We may choose ${\cal
D}_{\beta_i}$ so that its components are pairwise disjoint.  Note
however, that a component of ${\cal D}_{\beta_i}$ may have to
intersect $Q_{[s_i, f_i]}$ in its interior.

\begin{figure}
\centerline{\epsfxsize=3.0in \epsfbox{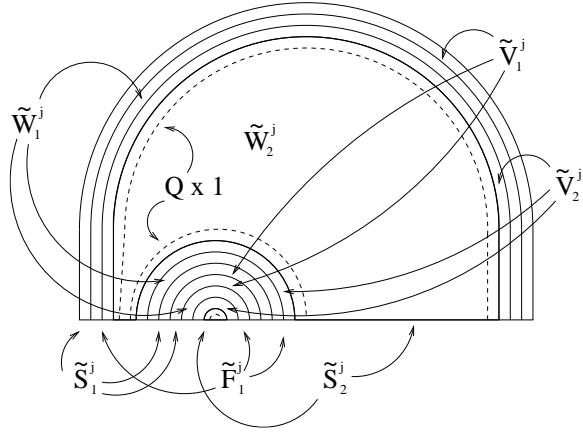}}
\caption{A boundary compressing disk intersecting $Q_{[s_2, f_2]}
\times I \subset Q \times I$}
\label{betadisk}
\end{figure}

\begin{figure}
\centerline{\epsfxsize=3.0in \epsfbox{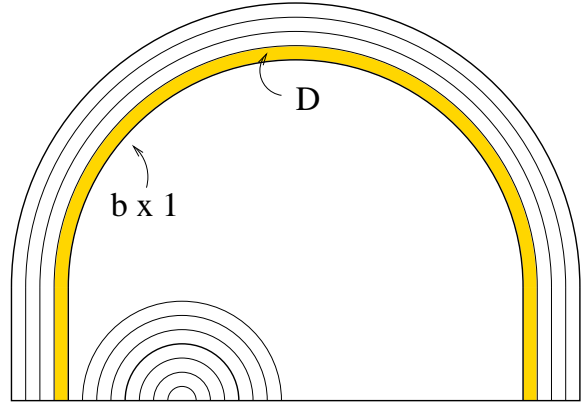}}
\caption{A destabilizing disk for the component $b$ of $\beta_2$}
\label{betadisk1}
\end{figure}

Let $D$ be a component of ${\cal D}_{\beta_i}$ and consider how $D$
meets $\tilde {\cal F} \cup \tilde {\cal S}$.  See Figure
\ref{betadisk} for the case $i = 2$.  Then the corresponding arc, call
it $b \times 1$, in $\beta _i$ is parallel to $\tilde S_i^j$ via the
disk sketched in Figure \ref{betadisk1}.  Extending $D$ into the other
side of $Q_{[s_i, f_i]}$ allows us to locate another such disk for $b
\times -1$.

The cocore of the 1-handle attached along the arc and a truncated
version of this disk define a destabilizing pair.  In this way, each
component of $\beta_i$ is in fact superfluous.
\qed

\vspace{1 mm}

The procedures above show that the number of components of
$Q_{[f_{i-1}, s_i]}$ that do not meet $q_{f_{i-1}}$ plays a role in
the complexities of the generalized Heegaard splittings constructed.
It is of particular importance to control the contribution arising from
annular components of this type.

\begin{defn}
An annulus $A$ in a compression body $W$ is called a
\underline{dipping annulus} if it is essential and $\partial A \subset
\partial_+W$.
\end{defn}

\begin{lem} \label{lem:dip}
Suppose ${\cal B}$ is a collection of essential annuli in a
compression body $W$ and that ${\cal A}$ is the subcollection
consisting of dipping annuli.  Denote the number of annular components
of $\partial_+W - {\cal B}$ by $l$.  If no two components of ${\cal
A}$ are isotopic, then \[|{\cal A}| \leq J(W) + \frac{l}{2}\]
\end{lem}

\proof
This is \cite[Lemma 7.3]{SS3}.  
\qed

\begin{cor} \label{cor:dip}
Suppose ${\cal B}$ is a collection of essential surfaces in a
compression body $W$ and that ${\cal A}$ is the subcollection
consisting of dipping annuli.  Denote the number of annular components
of $\partial_+W - {\cal B}$ that meet dipping annuli by
$l$.  Then \[|{\cal A}| \leq J(W) + \frac{l}{2}\]
\end{cor}

\proof First note that the conclusion depends only on ${\cal A}$, so
we may ignore ${\cal B - A}$ in the proof.  The proof is by induction
on the number of parallel dipping annuli.  If this number is $0$, then
the conclusion follows from Lemma \ref{lem:dip}.  To verify the
inductive step, observe that deleting a dipping annulus that is
parallel to another decreases the number of annular components of
$\partial_+W - {\cal B}$ by two.  \qed

\begin{figure}
\centerline{\epsfxsize=3.0in \epsfbox{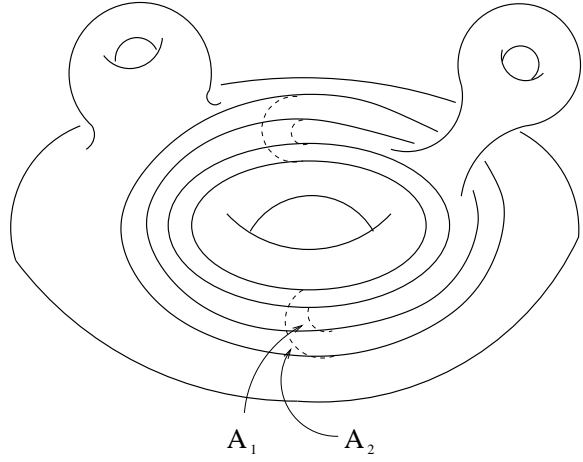}}
\caption{Isotopic dipping annuli}
\label{H12}
\end{figure}

\begin{defn}
Denote by $n_j$ the maximal number of pairwise non isotopic annuli
that can be simultaneously embedded in $M^j$.
\end{defn}

Recall also that $d_i$ is the number of components of $Q_{[f_{i-1},
s_i]}$ that do not meet $F_{i-1}$.

\begin{lem} \label{lem:gr}
In the counting arguments below, we may assume that in Lemma
\ref{lem:arcs1}
\[\sum_i d_i \leq \sum_i (-\chi(Q_{[f_{i-1}, s_i]}) + J(V_i)) + 2\sum_j
n_j\]
\end{lem}

\proof Recall that $d_i$ is the number of components of $Q_{[f_{i-1},
s_i]}$ that do not meet $\partial_-V_i$.  Denote the number of
components of $Q_{[f_{i-1}, s_i]}$ that are dipping annuli by $a_i$.
Then \[d_i \leq a_i - \chi(Q_{[f_{i-1}, s_i]})\] 

Corollary \ref{cor:dip} provides a bound on $a_i$, but this bound
depends on the number of annular components of $S_i - Q$ that meet
dipping annuli.  See for instance Figure \ref{H12}.  This number is
potentially unbounded, but we will show below that corresponding to
each annular component in $S_i - Q$ there is a destabilization of
$\dot V_1^j \cup_{\dot S_1^j} \dot W_1^j$.

\begin{figure}
\centerline{\epsfxsize=2.5in \epsfbox{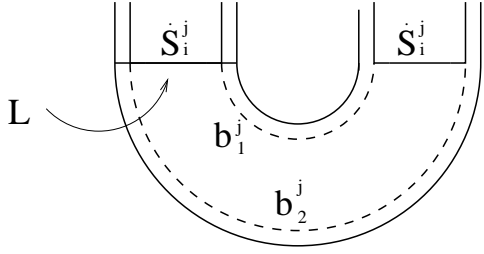}}
\caption{The annulus $L$}
\label{annulus}
\end{figure}

So consider an annulus $L$ in $\cup_i S_i - Q$ that meets at least one
dipping annulus.  See Figure \ref{annulus}.

\vspace{2 mm}
\noindent
Case I: $L$ is isotopic into $\partial M^j$. 

\begin{figure}
\centerline{\epsfxsize=2.5in \epsfbox{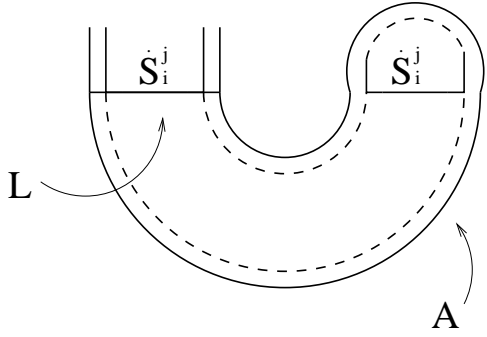}}
\caption{The annuli $L$ and $A$}
\label{bannulus}
\end{figure}

In this case there is an annulus $A$ in $\partial M^j$ as pictured
in Figure \ref{bannulus}. 

\begin{figure}
\centerline{\epsfxsize=3.5in \epsfbox{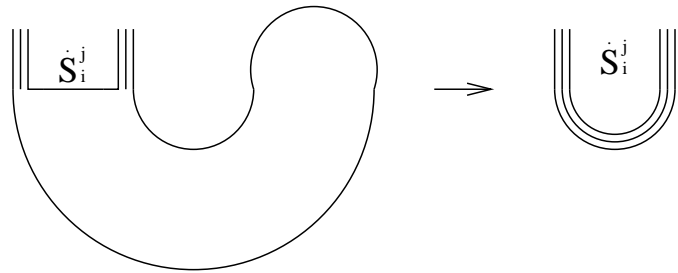}}
\caption{A simpler Heegaard splitting}
\label{trunc}
\end{figure}

Since $L$ is boundary parallel, $M^j$ is in fact homeomorphic to one
of the two components $C_1$ or $C_2$, say $C_1$, obtained by taking
the completion with respect to the path metric of $M^j - L$.  There is
thus a simpler generalized Heegaard splitting for $M^j$ than the one
currently under consideration.  This is the generalized Heegaard
splitting obtained from the one under consideration by deleting all
components of $\cup_i (\dot S_i^j \cup \dot F_i^j)$ in $C_2$ and
capping off any components of $\cup_i (\dot S_i^j \cup \dot F_i^j)$
with annuli.  See Figure \ref{trunc}.  The upshot is that
corresponding to each dipping annulus adjacent to $L$ there is a
destabilization of $\dot V_1^j \cup_{\dot S_1^j} \dot W_1^j$.  Note
that this is also true for any other annular components of $\cup_i S_i
- Q$ that meet dipping annuli and are contained in $C_2$.

\vspace{2 mm}
\noindent
Case II: $L$ is not isotopic into $\partial M^j$.

\begin{figure}
\centerline{\epsfxsize=3in \epsfbox{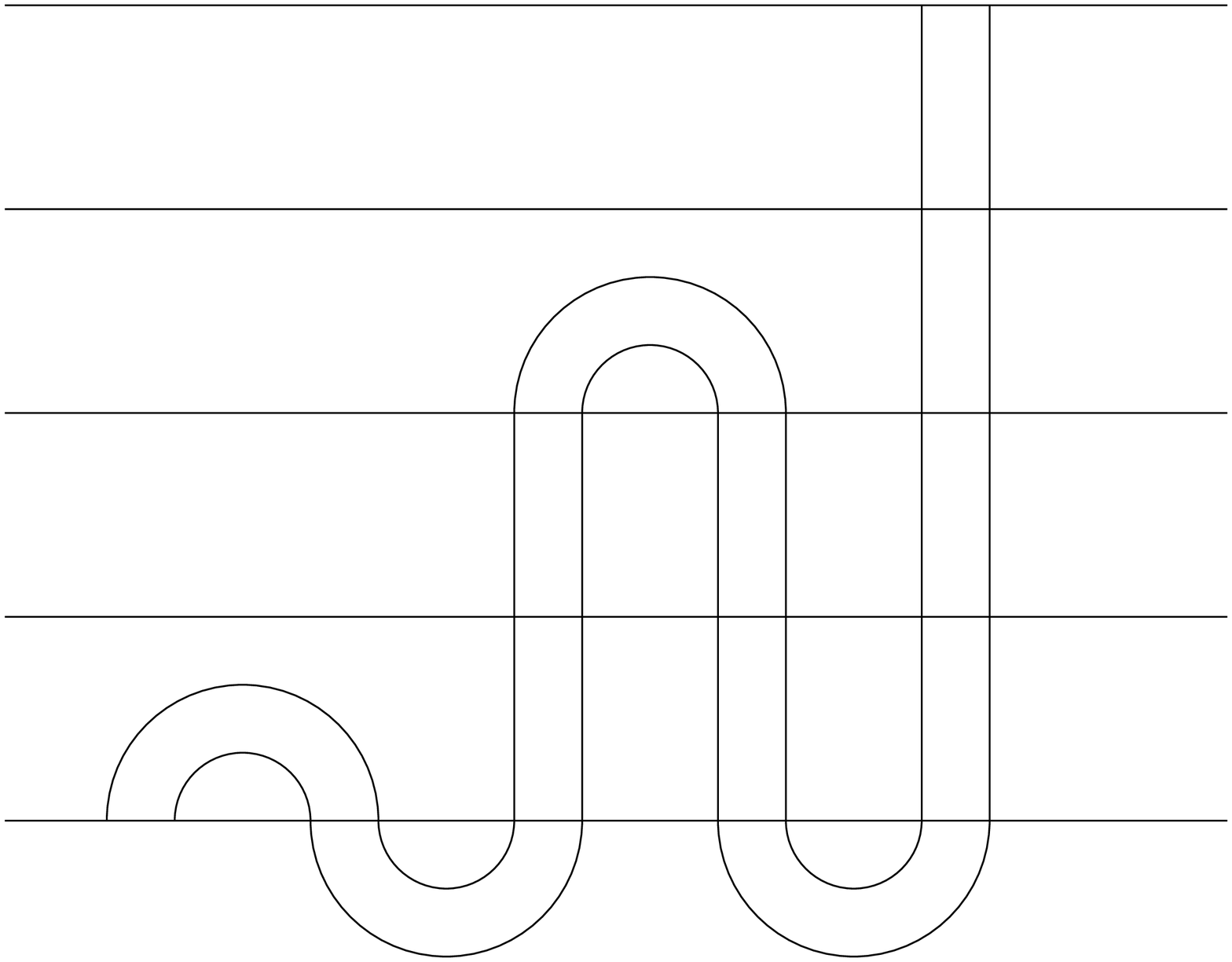}}
\caption{Schematic for $max(L)$}
\label{max(L)}
\end{figure}

\begin{figure}
\centerline{\epsfxsize=3in \epsfbox{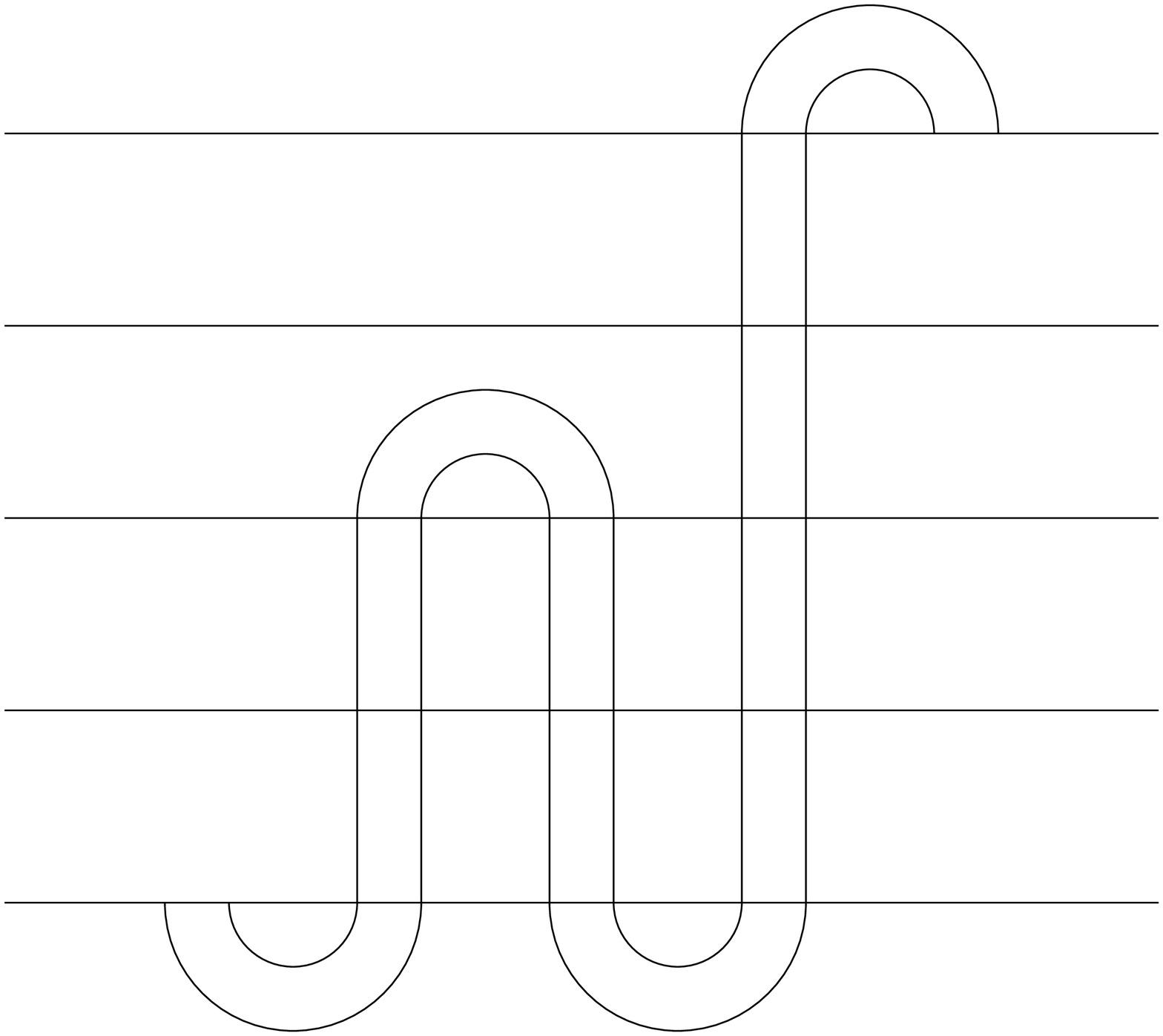}}
\caption{Schematic for $max(L)$}
\label{max(L)a}
\end{figure}

\begin{figure}
\centerline{\epsfxsize=1in \epsfbox{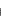}}
\caption{Schematic for a simpler generalized Heegaard splitting near $max(L)$}
\label{max(L)s}
\end{figure}

\begin{figure}
\centerline{\epsfxsize=1in \epsfbox{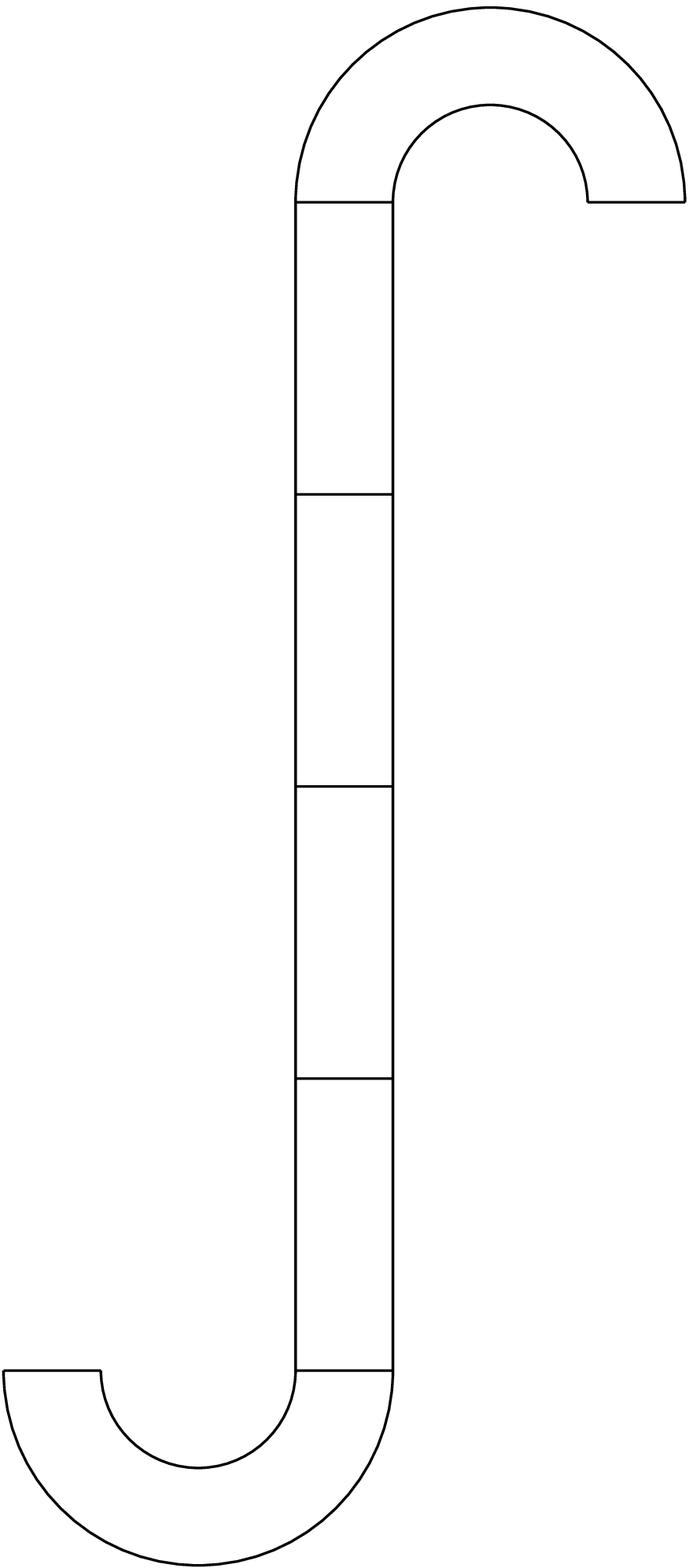}}
\caption{Schematic for a simpler generalized Heegaard splitting near $max(L)$}
\label{max(L)as}
\end{figure}

Denote by $max(L)$ the maximal product neighborhood of $L$ in $M - Q$
that is bounded by annular components of $\cup_i S_i - Q$.  See
Figures \ref{max(L)} and \ref{max(L)a}.  Then a simpler Heegaard
splitting may be constructed by replacing the portion of $\dot V_1^j
\cup_{\dot S_1^j} \dot W_1^j$ in $max(L)$ by one corresponding to a
simpler schematic.  See Figures \ref{max(L)s} and \ref{max(L)as}.

\begin{figure}
\centerline{\epsfxsize=1in \epsfbox{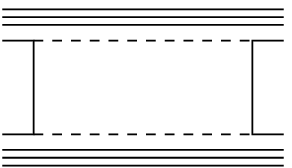}}
\caption{Before the ``straightening''}
\label{max(L)b}
\end{figure}

\begin{figure}
\centerline{\epsfxsize=1in \epsfbox{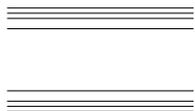}}
\caption{After the ``straightening''}
\label{max(L)ba}
\end{figure}

The effect of this ``straightening'' of $\dot V_1^j \cup_{\dot S_1^j}
\dot W_1^j$ in $max(L)$ near a pair of parallel dipping annuli
abutting $max(L)$ is pictured in Figures \ref{max(L)b} and \ref{max(L)ba}.

The upshot is that with the exception of at most one such pair of
 dipping annuli, each annular component of $\cup_i S_i - Q$ abutting
 dipping annuli corresponds to a destabilization.  
 By Corollary \ref{cor:dip}, this means that we may assume that
 
\[|{\cal A}| \leq \sum_i J(V_i) + 2\sum_j n_j\]

And hence that

\[\sum_i d_i \leq \sum_i (-\chi(Q_{[f_{i-1}, s_i]}) + a_i) \leq 
\sum_i (-\chi(Q_{[f_{i-1}, s_i]}) + J(V_i)) + 2\sum_j
n_j\]

 \qed

\section{Putting it all together}

By performing the construction, counting indices, subtracting amounts
corresponding to performing destabilizations and genus reductions, we
arrive at two propositions that imply the Main Theorem.

\begin{prop} \label{prop:par}
Let $M$ be a compact possibly closed orientable irreducible
$3$-manifold.  Let $M = (V_1 \cup_{S_1} W_1) \cup_{F_1} \dots
\cup_{F_{n-1}} (V_n \cup_{S_n} W_n)$ be a strongly irreducible
generalized Heegaard splitting.  Let $Q$ be a closed essential surface
isotopic to a subsurface of $\cal F$.  Then for $M^1, \dots, M^j$ the
completions of the components of $M - Q$ with respect to the path
metric there are generalized Heegaard splittings $M^j = (A_1^j
\cup_{G_1^j} B_1^j) \cup_{P_1^j} \dots \cup_{P_{n-1}^j} (A_n^j
\cup_{G_n^j} B_n^j)$ for $j = 1, \dots, k$ such that
\[\sum_i J(V_i) = \sum_j \sum_i J(A_i^j)\]
\end{prop}

\begin{figure}
\centerline{\epsfxsize=3.0in \epsfbox{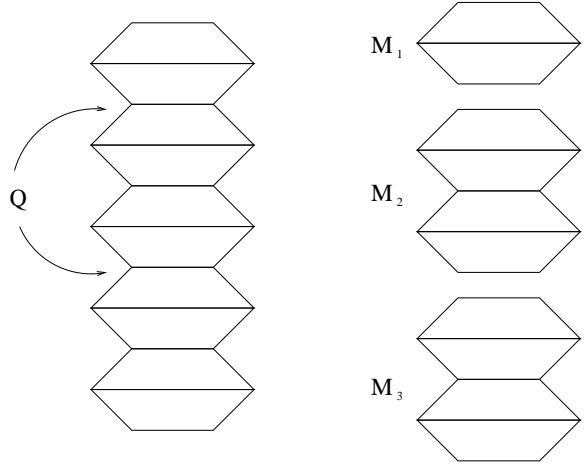}}
\caption{Cutting along a subset of ${\cal F}$}
\label{par}
\end{figure}

\proof
Isotope $Q$ to coincide with a subset of $\cal F$.  Then let $M^j$ be
the completion of a component of $M - Q$ and let $V_{i_1}, \dots,
V_{i_l}$, $W_{i_1}, \dots, W_{i_l}$ be the compression bodies among
$V_1, \dots, V_n, W_1, \dots, W_n$ that meet $M^j$.  Set $A_i^j =
V_i$, if $i \epsilon \{i_1, \dots, i_l\}$ and $A_i^j = \emptyset$
otherwise, $B_i^j = W_i$, if $i \epsilon \{i_1, \dots, i_l\}$ and
$B_i^j = \emptyset$ otherwise, $G_i^j = S_i$ if $i \epsilon \{i_1,
\dots, i_l\}$ and $G_i^j = \emptyset$ otherwise, and $P_i^j = F_i$ if
$i \epsilon \{i_1, \dots, i_l\}$ and $P_i^j = \emptyset$ otherwise.
See Figure \ref{par}.
\qed

\begin{prop} \label{prop:ess}
Let $M$ be a compact orientable irreducible
$3$-manifold.  Let $M = (V_1 \cup_{S_1} W_1) \cup_{F_1} \dots
\cup_{F_{n-1}} (V_n \cup_{S_n} W_n)$ be a strongly irreducible
generalized Heegaard splitting.  Let $Q$ be a compact
boundary incompressible essential surface in $M$.

Suppose that $\partial Q \subset \partial_-V$.  Suppose further that
no component of $Q$ is parallel to a component of $\cal F$.  Denote
the completions of the components of $M - Q$ with respect to the path
metric by $M^1, \dots, M^k$.  Further denote the number of pairwise
non isotopic annuli that can be embedded simultaneously in $M^j$ by
$n_j$.

Then there are generalized Heegaard splittings $M^j = (\dot V_1^j
\cup_{\dot S_1^j} \dot W_1^j) \cup_{\dot F_1^j} \dots \cup_{\dot
F_{n-1}^j} (\dot V_n^j \cup_{\dot S_n^j} \dot W_n^j)$ for $M^j$ for $j
= 1, \dots, k$ that satisfy the following inequality:
\[\sum_i J(V_i) \geq \frac{1}{5}(\sum_j \sum_i J(\dot V_i^j) + 6 \chi(Q) - 
8 \sum_j n_j)\]

\end{prop}

\proof Since $M = (V_1 \cup_{S_1} W_1) \cup_{F_1} \dots \cup_{F_{n-1}}
(V_n \cup_{S_n} W_n)$ is strongly irreducible $Q$ may be isotoped so
that each component of $Q \cap (\cal F \cup \cal S)$ is essential in
both $Q$ and $\cal F \cup \cal S$ and so that the number of components
of $Q \cap (\cal F \cup \cal S)$ is minimal.  By Lemma
\ref{lem:construction}, the Main Construction gives generalized
Heegaard splittings $M^j = (\dot V_1^j \cup_{\dot S_1^j}\dot W_1^j)
\cup_{\dot F_1^j} \dots \cup_{\dot F_{n-1}^j} (\dot V_n^j \cup_{\dot
S_n^j} \dot W_n^j)$ for which

\[\sum_j \sum_i J(\dot
V_i^j) = \sum_i (J(V_i) + 2\chi(Q_{[f_{i-1}, s_i]}) + 4(|\alpha_i| +
|\beta_i|))\]

By Lemma \ref{lem:arcs1} we may choose $\alpha_i$ so that 
\[\sum_i |\alpha_i| = \sum_i (-\chi(Q_{[f_{i-1}, s_i]}) + d_i)\] 

and by Lemma \ref{lem:arcs2} we may choose $\beta_i = \emptyset$.  Thus
\[\sum_i |\beta_i| = 0\]

By Lemma \ref{lem:gr}, we may assume that
\[\sum_i d_i \leq \sum_i (-\chi(Q_{[f_{i-1}, s_i]}) + J(V_i)) + 2\sum_j n_j\]

Thus \[\sum_i 4(|\alpha_i| + |\beta_i|) \leq 4 \sum_i
(-\chi(Q_{[f_{i-1}, s_i]}) + d_i) \leq \]
\[4 \sum_i (-2\chi(Q_{[f_{i-1}, s_i]}) + J(V_i)) + 8 \sum_j n_j\]

Thus 
\[\sum_j \sum_i J(\dot
V_i^j) \leq \sum_i (J(V_i) + 2\chi(Q_{[f_{i-1}, s_i]}) 
-8\chi(Q_{[f_{i-1}, s_i]}) + 4J(V_i)) + 8\sum_j n_j = \]
\[\sum_i (5J(V_i) - 6\chi(Q_{[f_{i-1}, s_i]})) + 8 \sum_j n_j\]

Whence \[\sum_i J(V_i) \geq \frac{1}{5}(\sum_j \sum_i J(\dot V_i^j) + 
6 \chi(Q) - 8 \sum_j n_j)\]

\qed

\begin{thm} \label{thm:main}
(The Main Theorem) Let $M$ be a compact orientable irreducible
$3$-manifold.  Let $M = (V_1 \cup_{S_1} W_1) \cup_{F_1} \dots
\cup_{F_{n-1}} (V_n \cup_{S_n} W_n)$ be a strongly irreducible
generalized Heegaard splitting.  Let $Q$ be a compact
boundary incompressible essential surface in $M$.

Suppose that $\partial Q \subset \partial_-V$.  Denote the completions
of the components of $M - Q$ with respect to the path metric by $M^1,
\dots, M^k$.  Further denote the number of pairwise non isotopic
annuli that can be embedded simultaneously in $M^j$ by $n_j$.

Then there are generalized Heegaard splittings $M^j = (A_1^j
\cup_{G_1^j} B_1^j) \cup_{P_1^j} \dots \cup_{P_{n-1}^j} (A_n^j
\cup_{G_n^j} B_n^j)$ for $M^j$ for $j = 1, \dots, k$ that satisfy the
following inequality
\[\sum_i J(V_i) \geq \frac{1}{5}(\sum_j \sum_i J(\dot A_i^j) + 
6 \chi(Q) - 8 \sum_j n_j)\]

\end{thm}

\proof By Lemma \ref{lem:essential}, $Q$ may be isotoped so that each
component of $Q \cap (\cal F \cup \cal S)$ is essential in both $Q$
and $\cal F \cup \cal S$.  We may assume that the number of components
of $Q \cap (\cal F \cup \cal S)$ is minimal subject to this condition.

Partition $Q$ into $Q^p \sqcup Q^n$, where $Q^p$ consists of those
components of $Q$ that are parallel to a component of ${\cal F}$ and
$Q^n$ consists of those components of $Q$ that are not parallel to any
component of ${\cal F}$.  Then proceed first as in Proposition
\ref{prop:par} using $Q^p$ instead of all of $Q$.  This yields an
equality.  In each of the resulting 3-manifolds proceed as in
Proposition \ref{prop:ess} using the appropriate subset of $Q^n$
instead of all of $Q$.  This yields the required inequality.  \qed

\begin{thm} \label{thm:genus}
Let $M$ be a compact orientable irreducible
$3$-manifold.  Let $Q$ be a boundary incompressible essential surface
in $M$.  Denote the completions of the components of $M - Q$ with
respect to the path metric by $M^1, \dots, M^k$.  Further denote the
number of pairwise non isotopic annuli that can be embedded
simultaneously in $M^j$ by $n_j$.

Then 
\[g(M, \partial M) \geq \frac{1}{5}(\sum_j g(M^j, \partial M^j) - 
|M - Q| + 5 - 2\chi(\partial_-V) + 4\chi(Q) - 4\sum_j n_j)\] 

\end{thm}

\proof 
Let $M = V \cup_S W$ be a Heegaard splitting that realizes
$g(M, \partial M)$.  Let $M = (V_1 \cup_{S_1} W_1) \cup_{F_1} \dots
\cup_{F_{n-1}} (V_n \cup_{S_n} W_n)$ be a weak reduction of $M = V
\cup_S W$.  Then 
\[\sum_{i=1}^m J(V_i) = 2g - 2 + \chi(\partial_-V)\]

Now apply Theorem \ref{thm:main} to $M = (V_1 \cup_{S_1} W_1)
\cup_{F_1} \dots \cup_{F_{n-1}} (V_n \cup_{S_n} W_n)$.  This yields
generalized Heegaard splittings $M^j = (A_1^j \cup_{G_1^j} B_1^j)
\cup_{P_1^j} \dots \cup_{P_{n-1}^j} (A_n^j \cup_{G_n^j} B_n^j)$ for
which

\[\sum_i J(V_i) \geq \frac{1}{5}(\sum_j \sum_i J(\dot A_i^j) + 
6 \chi(Q) - 8 \sum_j n_j)\]

Amalgamating $M^j = (A_1^j \cup_{G_1^j} B_1^j) \cup_{P_1^j} \dots
\cup_{P_{n-1}^j} (A_n^j \cup_{G_n^j} B_n^j)$ yields Heegaard
splittings $M^j = V^j \cup_{S^j} W^j$ with

\[\sum_i J(A_i^j) = 2g(S^j) - 2 + \chi(\partial_-V^j)\]

Hence

\[2g(S) - 2 + \chi(\partial_-V) = \sum_i J(V_i) \geq \]

\[\frac{1}{5}(\sum_j \sum_i J(\dot A_i^j) + 6 \chi(Q) - 8 \sum_j n_j) = \]

\[\frac{1}{5}(\sum_j (2g(S^j) - 2 + \chi(\partial_-V^j)) + 6 \chi(Q) - 8 \sum_j  n_j)\]

Here

\[\sum_j \chi(\partial_-V^j) = \chi(\partial_-V) + 2\chi(Q)\]

Therefore 

\[2g(S) - 2 + \chi(\partial_-V) \geq \frac{1}{5}(\sum_j (2g(S^j) - 2) + 
\chi(\partial_-V) + 8\chi(Q) - 8\sum_j n_j)\]

Hence 

\[2g(S) \geq \frac{1}{5}(\sum_j (2g(S^j) - 2) + 10 - 4\chi(\partial_-V) + 
8\chi(Q) -  8\sum_j n_j)\]

Whence

\[g(M, \partial M) \geq \frac{1}{5}(\sum_j g(M^j, \partial M^j) - 
|M - Q| + 5 - 2\chi(\partial_-V) + 4\chi(Q) - 4\sum_j n_j)\] 

\qed

Finally, we consider two interesting cases encompassed by this
construction.  In \cite{EM}, M. Eudave-Mu\~{n}oz constructs a family
of tunnel number $1$ knots whose complements contain incompressible
surfaces of arbitrarily high genus.  For elementary definitions
pertaining to knot theory, see for instance \cite{BZ1} or \cite{R}.

\begin{rem}
Let $K$ be a tunnel number one knot whose complement $C(K) = S^3 -
\eta(K)$ contains a closed incompressible surface $F$ of genus $g$.
Let $M^1, M^2$ be the manifolds obtained by cutting $C(K)$ along $F$.
(Here $M^1$, say, is the complement of an open regular neighborhood of
the toroidal graph in the construction of Eudave-Mu\~{n}oz and $M^2$
is the component containing $\partial C(K)$.)  In the examples of
Eudave-Mu\~{n}oz the genus is $g(F) + 1$ for both $M^1$ and $M^2$,
i.e., \[g(M^1) + g(M^2) = 2g(F) + 2\] This number is significantly
lower than the upper bound derived here.  This leads one to believe
that the bound suggested is too high.  However, the constructions in
this paper show that the upper bound in the inequality would be
obtained if the splitting surface of a minimal genus Heegaard
splitting had to intersect the incompressible surface many times.  The
examples of Eudave-Mu\~{n}oz are very well behaved in this respect, as
the splitting surface of the genus $2$ Heegaard splitting may be
isotoped to intersect the incompressible surface in only one curve
inessential in $F$ or two curves essential in $F$.
\end{rem}

\begin{rem}
Consider $R \times I$ for $R$ a closed connected orientable surface.
The genus of $R \times I$ is $g(R)$.  On the other hand, the genus of
the manifold $M$ obtained by identifying $R \times \{0\}$ to $R \times
\{1\}$ is $2$ for a suitably chosen gluing homeomorphism.  For details
on choosing such a gluing homeomorphism, see \cite{saul}.  Here the
right hand side in Theorem \ref{thm:genus} is negative unless $R$ is a
torus.  This is disconcerting, but merely means that the result is
automatically true.  Again, the genus $2$ Heegaard splitting in this
example is well behaved with respect to its intersection with $R$; it
too may be isotoped so that it intersects $R$ either in one curve
inessential in $R$ or two curves essential in $R$.
\end{rem}

Recall that in the appendix to \cite{SS3}, A. Casson provided a class
of $3$-manifolds each having a strongly irreducible generalized
Heegaard splitting and each containing an essential annulus such that
the number of dipping annuli can be arbitrarily large.  These examples
provide one reason to believe the upper bound provided here is not
larger than necessary.

\noindent
Department of Mathematics \\
1 Shields Avenue \\
University of California, Davis \\
Davis, CA 95616

\end{document}